\documentclass[3p,number,sort&compress]{elsarticle}
\journal{Nonlinear Analysis}
\usepackage{amsmath,amssymb,mathtools}
\usepackage{nicefrac,bm,upgreek,pdf14}
\usepackage[mathscr]{eucal}
\usepackage{hyperref}
\makeatletter
\def\ps@pprintTitle{%
   \let\@oddhead\@empty
   \let\@evenhead\@empty
   \def\@oddfoot{\reset@font\hfil\thepage\hfil}
   \let\@evenfoot\@oddfoot
}
\makeatother
\hypersetup{
  colorlinks=true,
  citecolor=dblue,
  linkcolor=dblue,
  frenchlinks=false,
  pdfborder={0 0 0},
  naturalnames=false,
  hypertexnames=false,
  breaklinks}
\newtheorem{theorem}{Theorem}
\newtheorem{lemma}[theorem]{Lemma}
\newtheorem{corollary}[theorem]{Corollary}
\newdefinition{remark}{Remark}
\newdefinition{example}{Example}
\newproof{pf}{Proof}
\newcommand{\df}{\coloneqq}
\newcommand{\RR}{\mathbb{R}}

\newcommand{\cK}{{K}}  
\newcommand{\cO}{{\mathscr{O}}}  
\newcommand{\abs}[1]{\lvert#1\rvert}

\newcommand{\norm}[1]{\lVert#1\rVert}

\newcommand{\dprec}{{\;\prec\!\!\prec\;}}

\begin{document}
\begin{frontmatter}

\title{A counterexample to a nonlinear version of the\\ Kre{\u\i}n--Rutman theorem
by R. Mahadevan}

\author{Ari Arapostathis}

\address{Department of Electrical and Computer Engineering,\\
The University of Texas at Austin, 2501 Speedway, EER~7.824,\\
Austin, TX~~78712, USA}
\ead{ari@ece.utexas.edu}

\begin{abstract}
In this short note we present a simple counterexample to
a nonlinear version of the Kre{\u\i}n--Rutman theorem reported
in [\emph{Nonlinear Anal.\ }\textbf{11} (2007), 3084--3090].
Correct versions of this theorem, and related results
for superadditive maps are also presented.
\end{abstract}

\begin{keyword}
Kre{\u\i}n--Rutman theorem \sep Positive operator
\sep Principal eigenvalue \sep Convex cone

\MSC primary 47H07 \sep 47H10 \sep 47J10; secondary 47B65
\end{keyword}
\end{frontmatter}

\section{Introduction}\label{Intro}

Kre{\u\i}n and Rutman in their seminal work \cite{KR-50} have studied
linear operators which leave invariant a cone in a Banach space.
Recall that an \emph{ordered Banach space} is a real Banach space $X$
with a cone $\cK$, a nontrivial closed subset of $X$ satisfying
\begin{enumerate}[(a)]
\item
$t\cK\subset\cK$ for all $t\ge0$, where $t\cK = \{tx \,\colon x\in\cK\}\,$;
\item
$\cK+\cK\subset\cK\,$;
\item
$\cK\cap(-\cK)=\{0\}$, where $-\cK = \{-x \,\colon  x\in\cK\}$.
\end{enumerate}
As usual, we write $x\preceq y$ if $y-x\in\cK$, and $x\prec y$ if
$x\preceq y$ and $x\ne y$.
When the interior of $\cK$, denoted as $\mathring{\cK}$, is nonempty, we call
$X$ a \emph{strongly ordered} Banach space.
We also write $x\dprec y$ if $y-x\in\mathring\cK$.
A continuous map $T\colon X\to X$ is
\begin{enumerate}[1.]
\item
\emph{positive} if $T(\cK)\subset\cK\,$;
\item
\emph{strictly positive} if $T(\cK\setminus\{0\})\subset\cK\setminus\{0\}\,$;
\item
\emph{strongly positive} if $T(\cK\setminus\{0\})\subset\mathring\cK\,$;
\item
\emph{order-preserving} or \emph{increasing} if $x\preceq y \implies T(x)  \preceq T(y)\,$;
\item
\emph{strictly order-preserving} if $x\prec y \implies T(x)  \prec T(y)\,$;
\item
\emph{strongly order-preserving} if $x\prec y \implies T(x)  \dprec T(y)\,$;
\item
\emph{homogeneous of degree one}, or \emph{1-homogeneous},
if $T(tx)=tT(x)$ for all $t\ge0$.
\end{enumerate}

The following nonlinear extension to the Kre{\u\i}n--Rutman theorem was reported
in \cite{Mahadevan-07}.
For a 1-homogeneous map $T\colon X\to X$
we say that $\lambda\in\RR$ is an eigenvalue of $T$ if there exists a nonzero
$x\in X$, such that $T(x) = \lambda x$.
Recall that a map $T \colon X \to X$ is called \emph{completely continuous} if
it is continuous and compact.

\begin{theorem}[{\cite[Theorem~2]{Mahadevan-07}}]\label{T-1}
Let $T \colon X \to X$ be an order-preserving, 1-homogeneous, completely
continuous map such that for some $u \in \cK$ and $M > 0$, $MT(u) \succeq u$.
Then there exist $\lambda > 0$ and $\Hat{x} \in \cK$, with
$\norm{\Hat{x}}=1$, such that
$T(\Hat{x}) = \lambda \Hat{x}$.
Moreover,
if $\mathring\cK\ne \varnothing$ and
$T$ is strongly positive and strictly order-preserving, the following hold.
\begin{enumerate}[(i)]
\item
$\Hat{x}$ is the unique unit eigenvector in $\cK\,$;
\item
$\lambda\ge |\lambda'|$ for any real eigenvalue $\lambda'$  of $T\,$;
\item
$\lambda$ is geometrically simple.
\end{enumerate}
\end{theorem}

It turns out that the assertions in (i) and (iii) are not true under the
assumptions of the theorem.
In Section~\ref{S2} we present a counterexample to the theorem in
\cite{Mahadevan-07} mentioned above.
In Section~\ref{S3} we review correct versions of this result.
With the exception of Section~\ref{S3.3} concerning superadditive maps,
the remaining results in Section~\ref{S3} are a combination of existing
results in the literature, and no originality is claimed.

I also wish to thank the anonymous referee who brought to my attention
the work in \cite{Chang-14}, which employs the notions of semi-strong
positivity and semi-strongly increasing maps, and improves upon the results in
\cite{Chang-09}.
It turns out that Theorem~\ref{T3.2.1} in Section~\ref{S3} is a variation of
Theorem~2.3 in \cite{Chang-14}.

\section{A Counterexample}\label{S2}
The following is an example of a strongly positive, strictly order-preserving,
1-homogeneous, continuous map $T$ on $\RR^{2}$ that has multiple positive
unit eigenvectors.

\begin{example}\label{EX1}
Let $\cK=\bigl\{x=(x_{1},x_{2})\in\RR^{2} \,\colon  x_{i}\ge0\,,~ i=1,2\bigr\}$.
Define
\begin{equation*}
\cK_{1}\,=\,\bigl\{x\in\cK \,\colon x_{1}>2 x_{2}\bigr\}\,,\qquad
\cK_{3}\,=\,\bigl\{x\in\cK \,\colon x_{2}>2 x_{1}\bigr\}\,,
\end{equation*}
and $\cK_{2}=\cK\setminus(\cK_{1}\cup\cK_{3})$.
Let
\begin{equation*}
T(x) \,\df\,\begin{cases}
\left(\begin{smallmatrix}2&2\\[3pt]1&1\end{smallmatrix}\right) x
& \text{if}~x\in\cK_{1}\,
\\[5pt]
3 x& \text{if}~x\in\cK_{2}\,\\[5pt]
\left(\begin{smallmatrix}1&1\\[3pt]2&2\end{smallmatrix}\right) x
&\text{if}~x\in\cK_{3}\,.
\end{cases}
\end{equation*}
It is clear that $T\colon \cK\to\cK$ is continuous, 1-homogeneous, and
strongly positive.
Also every element of $\cK_{2}$ is an eigenvector of $T$.

It remains to show that $T$ is strictly order-preserving.
We examine all possible cases:
\begin{enumerate}[(i)]
\item
If $x\,,y\in\cK_{i}$, $i=1,2,3$, and $x\prec y$, then it is clear that
$T(x)\prec T(y)$.
\item
Suppose $x\in\cK_{1}$, $y\in\cK_{3}$, and $x\prec y$.
Then we must have
\begin{equation}\label{E-e01}
2 x_{2} \,<\, x_{1}\,\le\, y_{1} \,<\, \frac{y_{2}}{2}\,.
\end{equation}
By \eqref{E-e01} we obtain that
\begin{equation}\label{E-e02}
x_{1}+x_{2}\,<\, \frac{3 x_{1}}{2} \,\le\, \frac{3 y_{1}}{2} 
\,<\,\frac{y_{1}+y_{2}}{2}\,.
\end{equation}
Since
$T(x)=\left(\begin{smallmatrix}2\\[3pt]1\end{smallmatrix}\right) (x_{1}+x_{2})$
and $T(y)=\left(\begin{smallmatrix}1\\[3pt]2\end{smallmatrix}\right) (y_{1}+y_{2})$,
it follows by \eqref{E-e02} that $T(x)\prec T(y)$.
Also,  if $x\succ y$, then $T(x)\succ T(y)$ by symmetry.
\item
Suppose $x\in\cK_{1}$, $y\in\cK_{2}$ and $x\prec y$.
Then we have
\begin{equation}\label{E-e03}
2 x_{2} \,<\, x_{1}\,\le\, y_{1} \,\le\, 2 y_{2}\,.
\end{equation}
It follows by \eqref{E-e03} that $2(x_{1}+x_{2}) < 3 y_{1}$
and $x_{1}+x_{2} < 3 y_{2}$.
Therefore $T(x)\prec T(y)$.
On the other hand, if $x\succ y$, then we have
\begin{equation}\label{E-e04}
x_{1} \;>\; 2 x_{2} \;\ge\; 2 y_{2}\;\ge\; y_{1} \,,
\end{equation}
and by \eqref{E-e04} we obtain
$2(x_{1}+x_{2}) > 3 y_{1}$ and $x_{1}+x_{2} > 3 y_{2}$.
Therefore $T(x)\succ T(y)$.
Also, by symmetry, if $x\in\cK_{3}$ and $y\in\cK_{2}$, then the strictly order-preserving
property holds. 
\end{enumerate}
It follows by (i)--(iii) that $T$ is strictly order-preserving.
\end{example}

\section{Existence and Uniqueness Results}\label{S3}
We denote by $\cK^{*}$ the dual cone, i.e.,
$\cK^{*} = \{x^*\in X^{*}\,\colon \langle x^*,x\rangle\ge 0\text{\ for all\ } x\in \cK\}$.
The dual cone $\cK^{*}$ might not satisfy $\cK^{*}\cap(-\cK^{*})=\{0\}$,
so is not necessarily a cone.
If $X$ is strongly ordered, then $x\in\mathring{\cK}$ if and only if
$\langle x^*,x\rangle>0$ for all nontrivial $x^*\in\cK^{*}$.

A cone $\cK\subset X$ is said to be \emph{generating} if $X=\cK-\cK$,
and it is said to be \emph{total} if $X$ equals the closure of $\cK-\cK$.
A strongly ordered Banach space is always generating.
A cone $\cK\subset X$  is called
\emph{normal} if there exists a positive constant $\gamma$
such that $\norm{x+y}\ge \gamma \norm{x}$ for all $x,y\in\cK$.

For a 1-homogeneous, continuous map $T\colon  X\to X$ we define,
as in \cite{Mallet-Paret-02,Ogiwara-95},
\begin{align*}
\norm{T}_{+} &\,\df\, \sup\;\bigl\{\norm{Tx} \,\colon x\in\cK\,,~
\norm{x}\le 1\bigr\}\,,
\\[5pt]
\Tilde{\varrho}_{+}(T)
&\,\df\, \lim_{n\to\infty}\; \norm{T^{n}}_{+}^{\nicefrac{1}{n}}\,,
\\[5pt]
\mu(x) &\,\df\, \limsup_{n\to\infty}\; \norm{T^{n}(x)}^{\nicefrac{1}{n}}\,,
\\[5pt]
\varrho_{+}(T) &\,\df\, \sup_{x\in\cK}\;\mu(x)\,,\\[5pt]
\Hat{r}(T) &\,\df\, \sup\;\bigl\{\lambda\ge0\,\colon \exists x\in K\setminus\{0\}
\text{\ with\ } T(x)=\lambda x\bigr\}\,.
\end{align*}
The quantities $\Tilde{\varrho}_{+}(T)$, $\varrho_{+}(T)$,
and $\Hat{r}(T)$, are referred to in \cite{Mallet-Paret-02} as
the \emph{Bonsall's cone spectral radius}, the \emph{cone spectral radius},
and the \emph{cone eigenvalue spectral radius} of $T$, respectively.
For a 1-homogeneous, continuous map $T\colon  X\to X$ we always have
\cite[Proposition~2.1]{Mallet-Paret-02}
\begin{equation*}
\Hat{r}(T)\,\le\,\varrho_{+}(T)\,\le\,\Tilde{\varrho}_{+}(T)\,<\,\infty\,.
\end{equation*}
Also, if $T$ is compact, then $\varrho_{+}(T)=\Tilde{\varrho}_{+}(T)$
\citep[Theorem~2.3]{Mallet-Paret-02}.
The equality $\varrho_{+}(T)=\Tilde{\varrho}_{+}(T)$ also holds in the absence
of compactness, provided that $T$ is order preserving and the cone $\cK$ is normal
\citep[Theorem~2.2]{Mallet-Paret-02}.

We summarize the main hypothesis:
\begin{description}
\item[(H1)]
$T \colon  X \to X$ is an order-preserving,
1-homogeneous, completely continuous map.
\end{description}

\subsection{Existence of an eigenvector in \texorpdfstring{$\cK$}{K}
with a positive eigenvalue}
\label{S3.1}

Existence of an eigenvector of $T$ in $\cK$ with a positive eigenvalue, i.e.,
the existence part of Theorem~\ref{T-1}, is asserted in
\cite[Theorem~3.1]{Cui-12} under the following weaker assumption.
\begin{description}
\item[(A1)]
There exist a non-zero $u=v-w$ with $v,w\in\cK$ and such that $-u\notin\cK$,
a positive constant $M$, and a positive integer $p$ such
that $MT^{p}(u) \succeq u$.
\end{description}
On the other hand, is a direct consequence of the more general result in
\cite[Theorem~2.1]{Nussbaum-81} that if $S \colon  X \to X$ satisfies (H1) and
\begin{description}
\item[(A2)]
The orbit $\cO(S,x) \df \{ S^{n}(x) \,\colon n=1,2,\dotsc\}$
of some $x\in\cK$ is unbounded,
\end{description}
then there exist a constant $t_{0}\ge1$ and $x_{0}\in\cK$, with $\norm{x_{0}}=1$,
satisfying $S(x_{0})= t_{0}x_{0}$.
It thus turns out that \cite[Theorem~3.1]{Cui-12}, and hence also
the existence part of \cite[Theorem~2]{Mahadevan-07}, are a direct consequence
of \cite[Theorem~2.1]{Nussbaum-81} and the following lemma.

\begin{lemma}\label{L-2}
Suppose $T \colon  X \to X$ satisfies \textup{(H1)} and \textup{(A1)}.
Let $\varepsilon>0$ be arbitrary, and define $S= (M+\varepsilon)^{\nicefrac{1}{p}}T$.
Then $\cO(S,v)$ is unbounded.
\end{lemma}

\begin{pf}
We argue by contradiction.
If $\cO(S,v)$ is bounded, then it is also relatively compact.
By the order-preserving property and $1$-homogeneity we obtain
$S^{kp}(v)\succeq S^{kp}(u)\succeq \bigl(1+\nicefrac{\varepsilon}{M}\bigr)^{k} u$.
Therefore any limit point $y$ of $\{S^{kp}(v) \,\colon k=1,2,\dotsc\}$
satisfies $y\succeq \bigl(1+\nicefrac{\varepsilon}{M}\bigr)^{k} u$
for all $k=1,2,\dotsc$, and since $-u\notin\cK$ this is not possible.
\qed
\end{pf}

It then follows by \cite[Theorem~2.1]{Nussbaum-81} and Lemma~\ref{L-2}
that, under Assumptions~(H1) and (A1), there exists
$x_{0}\in\cK$ with $\norm{x_{0}}=1$
and $\lambda_{0}\ge M^{-\nicefrac{1}{p}}$ such that $T(x_{0}) =\lambda_{0} x_{0}$.
It is also clear from the above discussion that, under (H1),
a necessary and sufficient
condition for the existence of a positive eigenvalue with eigenvector in $\cK$
is that $\Tilde{\varrho}_{+}(T)>0$.
We remark here, that the assumption that $X$ is strongly ordered,
$\cK$ is normal, and $\Tilde{\varrho}_{+}(T)>0$,
it is shown in \cite[Proposition~3.1.5]{Ogiwara-95} that
$\Tilde{\varrho}_{+}(T)=\Hat{r}(T)$.

\subsection{Uniqueness and simplicity of the positive eigenvalue}\label{S3.2}

We define
\begin{equation*}
\sigma_{+}(T) \,\df\, \bigl\{\lambda>0\,\colon
T(x)=\lambda x\,,~x\in \cK\setminus\{0\}\bigr\}\,.
\end{equation*}

Ogiwara introduced the property of \emph{indecomposability}
\cite[Hypothesis~A4]{Ogiwara-95} to obtain the following.
Suppose that $X$ is strongly ordered, $\cK$ is
normal, and $T$ satisfies (H1) and is indecomposable.
Then $\sigma_{+}(T)=\lbrace\lambda_0\rbrace$, i.e., a singleton,
$\lambda_0$ is a simple eigenvalue of $T$, and the corresponding eigenvector
lies in $\mathring{\cK}$
\cite[Theorem~3.1.1, Corollary~3.1.6]{Ogiwara-95}.

A significant improvement of the above result can be found in \cite{Chang-14}.
Chang defines \emph{semi-strong positivity} of $T$ in \cite[Definition~4.5]{Chang-09} by
\begin{equation*}
\forall\,x\in\partial{K}\setminus\{0\}\,,\ \exists\, x^*\in K^*
\text{\ such that\ } \bigl\langle x^*, T(x)\bigr\rangle\, >\, 0 \,=\,
\langle x^*,x\rangle\,.
\end{equation*}
Also $T$ is called \emph{semi-strongly increasing} in \cite[Definition~2.1]{Chang-14} if
\begin{equation*}
\forall\,x,y\in X\,,\text{\ with\ }
x-y\in\partial{K}\setminus\{0\}\,,\ \exists\, x^*\in K^*
\text{\ such that\ } 
\bigl\langle x^*, T(x)-T(y)\bigr\rangle\, >\, 0 \,=\,
\langle x^*,x-y\rangle\,.
\end{equation*}
Then normality of $\cK$ is relaxed to assert the following.
If $X$ is strongly ordered, and $T$ satisfies (H1) and is semi-strong positive, then
there exists a unique positive eigenvalue with eigenvector in $\mathring{\cK}$.
In addition, if $T$ is semi-strongly increasing, then the eigenvalue is simple
\cite[Theorem~2.3]{Chang-14}.
It is also shown that the indecomposability hypothesis of Ogiwara is
equivalent to the semi-strongly increasing property \cite[Theorem~4.3]{Chang-14}.

In the sequel, we only assume that $X$ is strongly ordered,
and that $T$ is 1-homogeneous and order preserving, and comment on
the uniqueness and simplicity of the eigenvalue, provided that
$\sigma_{+}(T)\ne \varnothing$.

Consider the following hypothesis:
\begin{description}
\item[(B1)]
If $x\in\partial{\cK}\setminus\{0\}$, then
$x-\beta T(x) \notin\cK$ for all $\beta>0$.
\end{description}
It is clear that semi-strong positivity implies (B1).
On the other hand, it
is straightforward to show that if two eigenvectors $x_{0}$ and $y_{0}$ lie in
$\mathring{\cK}$, then the associated eigenvalues are equal.
In turn, it is easy to show  that, under (B1), every eigenvector in $\cK$
with a positive eigenvalue
has to lie in $\mathring{\cK}$,
and, consequently, that the positive eigenvalue is unique.
Also, following for example the proof
in \cite[Lemma~3.1.2]{Ogiwara-95}, we can show that (B1) implies
$T(\mathring\cK)\subset\mathring\cK$.

Next, consider the hypothesis
\begin{description}
\item[(B2)]
If $x-y\in\partial{\cK}\setminus\{0\}$, then
$x-y-\beta\, \bigl(T(x)-T(y)\bigr) \notin\cK$ for all $\beta>0$.
\end{description}
Clearly, (B2) $\Rightarrow$ (B1).
Also (B2) is weaker than the strong order preserving property.
Under (B2), following the argument in the
proof of \cite[Theorem~3.1.1]{Ogiwara-95}, one readily shows that
if there exists a unit eigenvector in $\mathring\cK$, then it is unique.

We summarize the above assertions in the form of the following theorem.

\begin{theorem}\label{T3.2.1}
Let $X$ be strongly ordered, and $T\colon K\to K$ be an order-preserving,
1-homogeneous map with $\sigma_{+}(T)\ne \varnothing$.
\begin{enumerate}[(i)]
\item
If (B1) holds, then $T(\mathring\cK)\subset\mathring\cK$, $\sigma_{+}(T)$ is
a singleton, and all eigenvectors lie in $\mathring\cK$.
\item
If (B2) holds, then the unique eigenvalue in $\sigma_{+}(T)$ is simple.
\end{enumerate}
\end{theorem}

\subsection{Superadditive maps}\label{S3.3}
We say that $T\colon X\to X$ is \emph{superadditive} (\emph{superadditive on $\cK$})
if $T(x+y)\succeq T(x) + T(y)$ for all $x,y\in X$ ($x,y\in\cK$).
It is clear that a (strictly, strongly) positive superadditive map is
(strictly, strongly) order-preserving.

\begin{theorem}\label{T-3}
Let  $T\colon \cK\to\cK$ be a superadditive, 
1-homogeneous, completely continuous map such that $\Tilde{\varrho}_{+}(T)>0$.
Then there exists $\lambda_{0}>0$ and $x_{0}\in\cK$ with $\norm{x_{0}}=1$
such that $T(x_{0})=\lambda_{0} x_{0}$.
Moreover, if (B1) holds, then $x_{0}$ is the unique unit eigenvector of $T$ in $\cK$.
\end{theorem}

\begin{pf}
Existence follows from Section~\ref{S3.1}.
Suppose $x_{0}$ and $y_{0}$ are two distinct unit eigenvectors in $\cK$.
As mentioned in Section~\ref{S3.2},
hypothesis (B1) implies that $x_{0},y_{0}\in\mathring{\cK}$,
and therefore these eigenvectors have a common eigenvalue $\lambda_{0}>0$.
Hence there exists $\alpha>0$ such that
$x_{0}-\alpha y_{0}\in \partial{\cK}\setminus\{0\}$.
Since $T$ is superadditive, we obtain
\begin{equation*}
T(x_{0}-\alpha y_{0})\,\preceq\, T(x_{0})-\alpha T(y_{0})\,=\,
\lambda_0(x_{0}-\alpha y_{0})\,.
\end{equation*}
This contradicts (B1) unless $x_{0}-\alpha y_{0}=0$.
Uniqueness of a unit eigenvector in $\cK$ follows.
\qed
\end{pf}

\begin{remark}
For a superadditive map $T$, we have
\begin{equation*}
x-y -\beta\,T(x-y) \,\succeq\, x-y-\beta\,\bigl(T(x)-T(y)\bigr)
\end{equation*}
Therefore if $T$ satisfies (B1) it also satisfies (B2).
\end{remark}

Let $\cK_{+}\df\cK$, $\cK_{-}\df -\cK$, and define
\begin{equation*}
\sigma(T) \,\df\, \bigl\{\lambda\in\RR\,\colon
T(x)=\lambda x\,,~x\in X\setminus\{0\}\bigr\}\,.
\end{equation*}

\begin{corollary}
Let  $T \colon X \to X$ be a positive, superadditive, 
1-homogeneous, completely continuous map such that $\Tilde{\varrho}_{+}(T)>0$.
Assuming (B1), there exist unique unit eigenvectors $x_{+}\in\cK_{+}$
and $x_{-}\in\cK_{-}$ with positive eigenvalues $\lambda_{+}$ and
$\lambda_{-}$, respectively.
Moreover, $\lambda_{-}\ge\lambda_{+}$.
Also, if $\lambda\in\sigma(T)$, then $\abs{\lambda}\le\lambda_{+}$.
\end{corollary}

\begin{pf}
By Theorem~\ref{T-3}, $T$ has a unique eigenvector $x_{+}\in\cK_{+}$
corresponding to an eigenvalue $\lambda_{+}>0$, and moreover
$x_{+}\in\mathring\cK_{+}$.
Define $S(x) \df - T(-x)$.
By superadditivity $-T(-x)\succeq T(x)$, which implies that
$S(x_{+})\succeq T(x_{+})= \lambda_{+}x_{+}$
Therefore (A1) holds for $S$ which implies the existence of a unit eigenvector
$x_{-}\in\cK_{-}$ for $S$ with a positive eigenvalue $\lambda_{-}$.
Also property (B1) for $T$ implies that if
$x\in\partial{\cK}\setminus\{0\}$, then
\begin{equation*}
x-\beta S(x)\,\preceq\, x-\beta T(x) \notin\cK\,,
\end{equation*}
so that property (B1) also holds for $S$.
Thus uniqueness of $x_{-}$ follows by Section~\ref{S3.2}.
Let $\alpha>0$ be such that $x_{-}+\alpha x_{+}\in\partial{\cK}_{-}$.
By superadditivity,
\begin{equation*}
T(x_{-}+\alpha x_{+}) \,\succeq\, T(x_{-})+\alpha T(x_{+})\,=\,
\lambda_{-}x_{-}+\alpha\lambda_{+}x_{+}\,.
\end{equation*}
By the order-preserving property,
we have  $\lambda_{-}x_{-}+\alpha\lambda_{+}x_{+}\preceq 0$, which
implies that $\lambda_{-}\ge\lambda_{+}$.

Suppose $T(x)=\lambda x$ for some $x\in X\setminus(\cK_{+}\cup\cK_{-})$,
with $x\ne0$.
Let $\alpha>0$ be such that $x_{+}+\alpha x\in \partial\cK_{+}$.
Since $T^{2}$ is order preserving, we have
$\lambda_{+}^{2}x_{+} + \lambda^{2}\alpha x\succeq 0$, which is possible only if
$\lambda_{+}\ge\abs{\lambda}$. 
\qed
\end{pf}

We would also like to mention here the stability results reported
in \cite{ABK-JTP} concerning
strongly continuous semigroups of superadditive operators on Banach spaces.

\section*{Acknowledgements}
I wish to thank the anonymous referee whose comments helped to substantially
improve this paper.
This research was supported in part by the Army Research Office
through grant W911NF-17-1-001, and in part by the Office of Naval Research
through grant N00014-16-1-2956.


\end{document}